\documentclass{amsart}

\usepackage{amssymb}
\usepackage{amsmath}
\usepackage{amsthm}
\usepackage{amsfonts}

\usepackage{a4wide}
\usepackage{longtable}
\usepackage{multirow}
\usepackage{setspace}

\newtheorem{theorem}{Theorem}[section]

\theoremstyle{definition}
\newtheorem{definition}{Definition}[section]

\numberwithin{equation}{section}

\begin{document}


\title[Construction and classification of \(p\)-ring class fields]{Construction and classification of \(p\)-ring class fields \\ modulo \(p\)-admissible conductors}

\author{Daniel C. Mayer}
\address{Naglergasse 53\\8010 Graz\\Austria}
\email{algebraic.number.theory@algebra.at}
\urladdr{http://www.algebra.at}

\thanks{Research supported by the Austrian Science Fund (FWF): projects J0497-PHY and P26008-N25}



\date{Tuesday, December 29, 2020}


\begin{abstract}
Each \(p\)-ring class field \(K_f\)
modulo a \(p\)-admissible conductor \(f\)
over a quadratic base field \(K\)
with \(p\)-ring class rank \(\varrho_f\) mod \(f\)
is classified according to Galois cohomology and
\textbf{differential principal factorization type}
of all members of its associated heterogeneous multiplet
\(\mathbf{M}(K_f)=\lbrack(N_{c,i})_{1\le i\le m(c)}\rbrack_{c\mid f}\)
of dihedral fields \(N_{c,i}\)
with various conductors \(c\mid f\) having \(p\)-multiplicities \(m(c)\) over \(K\)
such that \(\sum_{c\mid f}\,m(c)=\frac{p^{\varrho_f}-1}{p-1}\).
The advanced viewpoint of \textbf{classifying} the entire collection \(\mathbf{M}(K_f)\),
instead of its individual members separately,
admits considerably deeper insight into the class field theoretic structure of ring class fields,
and the actual \textbf{construction} of the multiplet \(\mathbf{M}(K_f)\)
is enabled by exploiting the routines for abelian extensions
in the computational algebra system Magma.
\end{abstract}

\maketitle


\section{Introduction}
\label{s:Intro}
\noindent
The aim of this article is
to present an entirely new technique
for the construction and classification
of non-Galois fields \(L\) of odd prime degree \(p\)
as subfields \(L<K_f\) of a \(p\)-ring class field \(K_f\)
modulo a \(p\)-admissible conductor \(f\)
over a quadratic base field \(K\).
The innovative idea underlying this new method
is the fact that,
if the Galois closure \(N\) of such a field \(L\)
is absolutely dihedral of degree \(2p\)
with automorphism group
\(\mathrm{Gal}(N/\mathbb{Q})\simeq
D_p=\langle\sigma,\tau\mid\sigma^p=\tau^2=1,\tau\sigma=\sigma^{-1}\tau\rangle\),
then \(N\) is relatively cyclic of degree \(p\)
with group \(G=\mathrm{Gal}(N/K)\simeq C_p=\langle\sigma\rangle\)
over its unique quadratic subfield \(K=\mathrm{Fix}(\sigma)\)
and can be viewed as an \textit{abelian extension}
modulo some conductor \(f\) over \(K\)
within the scope of class field theory
\cite{Ar1927,Be2004,Fi2001,Ha1930}.

The construction process for the fields \(L\) is implemented
as a program script for the
computational algebra system Magma
\cite{BCP1997,BCFS2020,MAGMA2020}
using the \textit{class field theoretic routines} by Fieker
\cite{Fi2001},
and the normal fields \(N/L\) are classified
according to the cohomology \(\hat{H}^0(G,U_N)\) and \(H^1(G,U_N)\)
of their unit group \(U_N\)
as a Galois module over \(G\)
\cite{Ma2019a,Ma2019b,Mo1979}.

For \(p\ge 5\), the results are completely new,
whereas for \(p=3\), they admit an independent verification
and a class field theoretic illumination
of classical tables of cubic fields by
Angell 1972
\cite{Angell1972,Angell1973}
and 1975
\cite{Angell1975,Angell1976},
Ennola and Turunen 1983
\cite{EnTu1983,EnTu1985},
Llorente and Quer 1988
\cite{LlQu1988},
Fung and Williams 1990
\cite{FuWi1990},
and Belabas 1997
\cite{Be1997}.
However, in contrast to these well-known tables,
where the focus was on the computation of
fundamental systems of units and the structure of ideal class groups
\cite{Angell1972,Angell1973,Angell1975,Angell1976,EnTu1983,EnTu1985,FuWi1990},
or even only of generating polynomials and prime decompositions
\cite{Be1997,LlQu1988},
our innovative database establishes an arrangement
according to conductors with an increasing number of prime factors,
pays attention to the phenomenon of \textit{multiplicities of discriminants}
\cite{Ma1992a,Ma2001,Ma2012a,Ma2012b,Ma2014},
and constitutes the \textit{first classification into} \(9\), respectively \(3\),
\textit{differential principal factorization types}
of totally real, respectively simply real, cubic number fields
\cite{Ma1991b,Ma1991c,Ma2019a,Ma2019b,Mo1979}.
This is a progressive new kind of structural information
which has never been provided for algebraic number fields before,
except for pure cubic fields
\cite{AMITA2020,BaCo1971,Wi1982}
and pure quintic fields
\cite{Ma2019a},
but the present paper emphasizes the advanced viewpoint
of \textbf{classifying an entire ring class field} \(K_f\)
by its associated \textit{heterogeneous multiplet} \(\mathbf{M}(K_f)\)
of dihedral fields with various conductors \(c\mid f\).


\section{Heterogeneous multiplets of objects and invariants}
\label{s:HeteroMultiplets}
\noindent
Let \(K=\mathbb{Q}(\sqrt{d})\) be a quadratic base field
with positive or negative fundamental discriminant \(d=d_K\equiv 0,1\,(\mathrm{mod}\,4)\),
essentially squarefree except possibly for the \(2\)-contribution \(v_2(d)\).
Suppose that \(p\) is an odd prime number and
\(f\ge 1\) is a \(p\)-admissible conductor over \(K\)
\cite{Ma1992a,Ma2014}.
Then the \(p\)-ring class field \(K_{p,f}\) mod \(f\)  of \(K\)
contains all cyclic relative extensions \(N/K\) with some conductor \(c\mid f\)
which are absolutely dihedral with automorphism group \(\mathrm{Gal}(N/\mathbb{Q})\simeq D_p\)
over the rational number field \(\mathbb{Q}\).
The crucial concept underlying this entire paper
is the collection of all these dihedral fields in a \textit{heterogeneous multiplet}
\(\mathbf{M}(K_{p,f})=\lbrack(N_{c,i})_{1\le i\le m_p(K,c)}\rbrack_{c\mid f}\)
according to the \(p\)-multiplicities \(m_p(K,c)\)
\cite{Ma1992a,Ma2014},
which satisfy the relation \(\sum_{c\mid f}\,m_p(K,c)=\frac{p^{\varrho_{p,f}}-1}{p-1}\)
in terms of the \(p\)-ring class rank \(\varrho_{p,f}\) modulo \(f\) of \(K\).
Since our principal aim is the classification of \(p\)-ring class fields \(K_{p,f}\),
it is essential to distinguish between
a multiplet of \textit{objects}
(expressing the multiplicity of the discriminants \(d_N\))
and a corresponding multiplet of \textit{invariants}
(expressing the Galois cohomology of the unit groups \(U_N\)
and differential principal factorizations of the fields \(N\)).

\begin{definition}
\label{dfn:Type}
By the \textbf{type of the \(p\)-ring class field} \(K_{p,f}\) modulo \(f\) of \(K\)
we understand the pair \((\mathrm{Obj}(K_{p,f}),\mathrm{Inv}(K_{p,f}))\) of heterogeneous multiplets
\begin{equation}
\label{eqn:Type}
\begin{aligned}
\mathrm{Obj}(K_{p,f}) &= \lbrack(N_{c,i})_{1\le i\le m_p(K,c)}\rbrack_{c\mid f} \\
\mathrm{Inv}(K_{p,f}) &= \lbrack(\tau(N_{c,i}))_{1\le i\le m_p(K,c)}\rbrack_{c\mid f}
\end{aligned}
\end{equation}
consisting of all absolutely dihedral fields \(N_{c,i}\) with conductors \(c\) dividing \(f\) as \textit{objects}
and their differential principal factorization (DPF) types \(\tau(N_{c,i})\) as \textit{invariants}
\cite{Ma2019a,Ma2019b}.
\end{definition}


\section{Homogeneous multiplets of unramified extensions}
\label{s:HomoMultiplets}
\noindent
The unique situation
where the heterogeneous multiplets degenerate to \textit{homogeneous multiplets}
occurs for \textit{unramified} relative extensions \(N/K\)
with conductor \(f=1\) which has only itself as a divisor \(c\mid f\).
In this unramified case,
which implies positive \(p\)-class rank \(\varrho_p=\varrho_{p,1}\ge 1\) of the quadratic base field \(K\),
there occur \textit{at most two} possible differential principal factorization types.


\begin{theorem}
\label{thm:Unramified}
\noindent
An \textbf{unramified} cyclic extension \(N\) with odd prime degree \(p\) of \(K\)
possesses the conductor \(f=1\) without any prime divisors.
For a \textbf{totally real} field \(N\), there are two cases:
\begin{enumerate}
\item
If the \(p\)-class rank of \(K\) is \(\varrho_p=1\),
then \(N\) is of \textbf{type \(\delta_1\)}.
\item
If the \(p\)-class rank of \(K\) is \(\varrho_p\ge 2\),
then \textbf{two types \(\alpha_1\) and \(\delta_1\)} are possible for \(N\).
\end{enumerate}
If \(N\) is \textbf{totally complex}, then \(N\) is of \textbf{type \(\alpha_1\)},
independently of the \(p\)-class rank of \(K\).
\end{theorem}

\begin{proof}
Since the conductor \(f=q_1\cdots q_t\) is essentially
the square free product of all prime numbers \(q_i\in\mathbb{P}\),
whose overlying prime ideals \(\mathfrak{q}_i\in\mathbb{P}_K\) are ramified in \(N\),
the following chain of equivalent statements is true:
\(N/K\) is unramified
\(\Longleftrightarrow\) None of the prime ideals of \(K\) ramifies in \(N\)
\(\Longleftrightarrow\) The conductor \(f=1\) has no prime divisors, i.e., \(t=0\).

Now we use the fundamental equation in
\cite[Cor. 5.1]{Ma2019b}
and the estimates in
\cite[Cor. 5.2]{Ma2019b}
for the decision about possible types of principal factorizations.
If \(f=1\), then there neither exist absolute principal factorizations in \(L/\mathbb{Q}\),
since \(0\le A\le\min(t,2)=0\),
nor relative principal factorizations in \(N/K\),
since \(0\le R\le\min(s,2)=0\),
where \(s\le t\) denotes the number of prime divisors \(q_i\) of \(f\) which split in \(K\).
Consequently, the fundamental equation degenerates to \(U+1=C\) with \(1\le U+1\le 2\),
which implies \(1\le C\le\min(\varrho_p,2)\).
Thus, only type \(\delta_1\) with \(C=1\) is possible for \(\varrho_p=1\),
whereas type \(\alpha_1\) with \(C=2\) can arise additionally for \(\varrho_p\ge 2\).
\end{proof}


\section{Conductors with a single prime divisor}
\label{s:SinglePrime}
\noindent
For a \textit{regular prime} conductor \(f\),
only two cases are possible.


\begin{theorem}
\label{thm:SingleRegularPrime}
\noindent
Let \(K\) be a quadratic base field with \(p\)-class rank \(\varrho=\varrho_p\).
Suppose \(f=q\) is a \textbf{regular} \(p\)-admissible \textbf{prime conductor} for \(K\).
Then the heterogeneous multiplet \(\mathbf{M}(K_{p,f})\) associated with the
\(p\)-ring class field \(K_{p,f}\) mod \(f\) of \(K\)
consists of two homogeneous multiplets with multiplicities
\(m_p(K,1)\) and \(m_p(K,q)\).
In this order,
and in dependence on the \(p\)-ring space \(V_p(q)\),
these two multiplicities are given by
\begin{enumerate}
\item
\((1+p+\ldots+p^{\varrho-1},\ p^\varrho)\), \quad if \(V_p(q)=V\) (\textbf{free} situation),
\item
\((1+p+\ldots+p^{\varrho-1},\ 0)\), \quad if \(V_p(q)<V\) (\textbf{restrictive} situation).
\end{enumerate}
\end{theorem}

\begin{proof}
See
\cite[Thm. 3.2, p. 2215, and Thm. 3.3, p. 2217]{Ma2014}.
\end{proof}


In the special case \(p=3\),
there also exists the possibility of an \textit{irregular prime power} conductor \(f=3^2\),
provided the discriminant of the quadratic field satisfies the congruence
\(d\equiv -3\,(\mathrm{mod}\,9)\).

\begin{theorem}
\label{thm:SingleIrregularPrime}
\noindent
Assume that \(p=3\).
Let \(K\) be a quadratic base field with \(3\)-class rank \(\varrho=\varrho_3\)
and discriminant \(d\equiv -3\,(\mathrm{mod}\,9)\).
Consider the \textbf{irregular} \(3\)-admissible \textbf{prime power conductor} \(f=3^2\) for \(K\).
Then the heterogeneous multiplet \(\mathbf{M}(K_{p,f})\) associated with the
\(3\)-ring class field \(K_{3,f}\) mod \(f\) of \(K\)
consists of three homogeneous multiplets with multiplicities
\(m_3(K,1)\), \(m_3(K,3)\) and \(m_3(K,9)\).
In this order,
and in dependence on the \(3\)-ring spaces \(V_3(3)\) and \(V_3(9)\),
these three multiplicities are given by
\begin{enumerate}
\item
\((1+3+\ldots+3^{\varrho-1},\ 3^\varrho,\ 3^{\varrho+1})\), \quad if \(V_3(9)=V_3(3)=V\) (\textbf{free} situation),
\item
\((1+3+\ldots+3^{\varrho-1},\ 3^\varrho,\ 0)\), \quad if \(V_3(9)<V_3(3)=V\),
\item
\((1+3+\ldots+3^{\varrho-1},\ 0,\ 3^\varrho)\), \quad if \(V_3(9)=V_3(3)<V\),
\item
\((1+3+\ldots+3^{\varrho-1},\ 0,\ 0)\), \quad if \(V_3(9)<V_3(3)<V\) (\textbf{maximal restriction}).
\end{enumerate}
\end{theorem}

\begin{proof}
See
\cite[Thm. 3.4, p. 2217]{Ma2014}.
\end{proof}


\section{Conductors with two prime divisors}
\label{s:TwoPrimes}
\noindent
For \textit{regular} conductors \(f\) divisible by \textit{two primes},
more distinct situations may arise.


\begin{theorem}
\label{thm:TwoPrimes}
\noindent
Let \(K\) be a quadratic base field with \(p\)-class rank \(\varrho=\varrho_p\).
Suppose \(f=q_1\cdot q_2\) is a \textbf{regular} \(p\)-admissible conductor for \(K\)
with \textbf{two prime divisors} \(q_1\) and \(q_2\).
Then the heterogeneous multiplet \(\mathbf{M}(K_{p,f})\) associated with the
\(p\)-ring class field \(K_{p,f}\) mod \(f\) of \(K\)
consists of four homogeneous multiplets with multiplicities
\(m_p(K,1)\), \(m_p(K,q_1)\), \(m_p(K,q_2)\) and \(m_p(K,f)\).
In this order,
and in dependence on the \(p\)-ring spaces \(V_p(q_1)\), \(V_p(q_2)\) and \(V_p(f)\),
these four multiplicities are given by
\begin{enumerate}
\item
\((1+p+\ldots+p^{\varrho-1},\ p^\varrho,\ p^\varrho,\ p^\varrho(p-1))\), \quad if \(V_p(f)=V_p(q_1)=V_p(q_2)=V\) (\textbf{free} case),
\item
\((1+p+\ldots+p^{\varrho-1},\ p^\varrho,\ 0,\ 0)\), \quad if \(V_p(f)=V_p(q_2)<V_p(q_1)=V\),
\item
\((1+p+\ldots+p^{\varrho-1},\ 0,\ p^\varrho,\ 0)\), \quad if \(V_p(f)=V_p(q_1)<V_p(q_2)=V\),
\item
\((1+p+\ldots+p^{\varrho-1},\ 0,\ 0,\ p^\varrho)\), \quad if \(V_p(f)=V_p(q_1)=V_p(q_2)<V\),
\item
\((1+p+\ldots+p^{\varrho-1},\ 0,\ 0,\ 0)\), \quad if \(V_p(f)<V_p(q_1)\ne V_p(q_2)<V\) (\textbf{maximal restriction}).
\end{enumerate}
\end{theorem}

\begin{proof}
We use the terminology and notation in
\cite{Ma2014}.
Generally, the \(p\)-ring class rank is given by
\(\varrho_{p,f}=\varrho+t+w-\delta_p(f)\).
Here, we have either \(t=2\), \(w=0\) or \(t=1\), \(w=1\),
and thus \(\varrho_{p,f}=\varrho+2-\delta_p(f)\).
Also, we know that generally
\(m_p(K,1)=\frac{p^\varrho-1}{p-1}\).
Since \(f=q_1\cdot q_2\) is \(p\)-admissible,
\(q_1\) and \(q_2\) must also be \(p\)-admissible, both.
\begin{enumerate}
\item
In the free case with defect \(\delta_p(f)=0\),
we have \(V_p(f)=V_p(q_1)=V_p(q_2)=V\) and
\[
\frac{p^{\varrho+2}-1}{p-1}-\frac{p^\varrho-1}{p-1}
=\frac{p^\varrho(p^2-1)}{p-1}
=p^\varrho(p+1)
=p^\varrho+p^\varrho+p^\varrho(p-1),
\]
which is exactly the desired partition
\[
\frac{p^{\varrho_{p,f}}-1}{p-1}-m_p(K,1)=m_p(K,q_1)+m_p(K,q_2)+m_p(K,f).
\]
\item
If \(q_1\) is free and \(q_2\), \(f\) are restrictive, then
\(V_p(f)=V_p(q_2)<V_p(q_1)=V\) and the relation
\[
\frac{p^{\varrho+1}-1}{p-1}-\frac{p^\varrho-1}{p-1}
=\frac{p^\varrho(p-1)}{p-1}
=p^\varrho,
\]
must be interpreted as \(m_p(K,q_1)=p^\varrho\) and \(m_p(K,q_2)=m_p(K,f)=0\).
\item
This case arises by interchanging the roles of \(q_1\) and \(q_2\) in the previous case.
\item
Additionally to (2) and (3),
there is another case of defect \(\delta_p(f)=1\)
where neither \(q_1\) nor \(q_2\) is free
but their \(p\)-ring spaces coincide \(V_p(f)=V_p(q_1)=V_p(q_2)<V\).
Then the formula in (2) has to be interpreted as 
\(m_p(K,q_1)=m_p(K,q_2)=0\) and \(m_p(K,f)=p^\varrho\).
\item
Finally, in the case of maximal restriction
with defect \(\delta_p(f)=2\),
which occurs for distinct \(p\)-ring spaces \(V_p(f)<V_p(q_1)\ne V_p(q_2)<V\),
there is no rank increment from \(\varrho\) to \(\varrho_{p,f}\),
and thus \(m_p(K,q_1)=m_p(K,q_2)=m_p(K,f)=0\).\qedhere
\end{enumerate}
\end{proof}


\section{Construction of \(p\)-ring class fields}
\label{s:Construction}
\noindent
This section describes
how the classification of non-trivial \(p\)-ring class fields
is prepared by their \textit{construction} and \textit{rigorous count}.
The intended class field theoretic illumination
of the structure of heterogeneous multiplets
\(\mathbf{M}(K_{p,f})=\lbrack(N_{c,1},\ldots,N_{c,m(c)})\rbrack_{c\mid f}\)
associated with \(p\)-ring class fields \(K_{p,f}\)
modulo \(p\)-admissible conductors \(f\)
over quadratic fields \(K\)
must pay \textit{primary attention} to the \(p\)-\textit{class rank} \(\varrho_p\)
of the quadratic base fields \(K=\mathbb{Q}(\sqrt{d})\),
since \(\varrho_p\) enters the formula for the multiplicities \(m(c)\).
More precisely,
since the existence of a torsion free fundamental unit \(\varepsilon>1\)
in real quadratic fields \(K\) with \(d>0\),
and the occurrence of the \(3\)-torsion unit \(\zeta_3\)
in the particular imaginary quadratic field \(K\) with \(d=-3\) in the case \(p=3\),
exerts a crucial impact on the codimension of \(p\)-ring spaces \(V_p(c)\),
the invariant \(\varrho_p\) must rather be replaced
by the \(p\)-\textit{Selmer rank} \(\sigma_p\) of \(K\)
which describes all \(p\)-virtual units of \(K\),
those which arise from non-trivial \(p\)-classes and the units in the usual sense:
\begin{equation}
\label{eqn:SelmerRank}
\sigma_p=
\begin{cases}
\varrho_p   & \text{ if } p\ge 5,\ d<0 \text{ or } p=3,\ d<-3, \\
\varrho_p+1 & \text{ if } d>0 \text{ or } p=3,\ d=-3.
\end{cases}
\end{equation}
The \textit{secondary attention} is devoted to various \(p\)-admissible conductors
\(f=q_1\cdots q_t\) with an increasing number \(t\ge 0\) of prime divisors,
starting with unramified extensions having \(t=0\), \(f=1\),
and continuing with ramified extensions,
beginning with prime or prime power conductors having \(t=1\), \(f=q_1\)
with a prime \(q_1\in\mathbb{P}\) or the critical prime power \(q_1=p^2\).


\section{Multiplets over imaginary quadratic fields for \(p=3\)}
\label{s:MultipletsComplex3}

\noindent
The focus of this section and most of the further sections is on \(p=3\),
where the components \(N_{c,i}\) of multiplets are cyclic cubic extensions
of quadratic base fields \(K\).
Here, we begin with imaginary base fields \(K\)
having the smallest possible \(3\)-Selmer rank \(\sigma_3=\varrho_3\).
The behavior of the particular imaginary quadratic field \(K\) with \(d=-3\)
where the extensions  \(N_{c,i}/K\) contain pure cubic fields
is rather similar to real quadratic base fields \(K\) with \(\sigma_3=\varrho_3+1\),
and thus the case \(d=-3\) will be treated separately.

\begin{theorem}
\label{thm:ComplexRank0}
Let \(K\) be an imaginary quadratic field
with fundamental discriminant \(d<-3\)
and trivial \(3\)-class rank \(\varrho_3=0\).
Assume that \(f=q_1\cdots q_{\tau}\) is a \(3\)-admissible conductor
with \(\tau\ge 1\) regular prime or prime power divisors \(q_i\)
(that is, either \(q_i\equiv\pm 1\,(\mathrm{mod}\,3)\)
or \(q_{\tau}=3\), \(d\equiv\pm 3\,(\mathrm{mod}\,9)\)
or \(q_{\tau}=9\), \(d\equiv\pm 1\,(\mathrm{mod}\,3)\)
but not \(q_{\tau}=9\), \(d\equiv -3\,(\mathrm{mod}\,9)\)).
Then the \(3\)-ring class field \(K_{3,f}\) modulo \(f\) of \(K\)
contains a homogeneous multiplet \(\mathbf{M}(K_{3,f})=(N_{f,1},\ldots,N_{f,m})\)
of dihedral fields with conductor \(f\) and multiplicity \(m=2^{\tau-1}\)
(singlet, doublet, quartet, octet, hexadecuplet, etc.).
\end{theorem}

\begin{proof}
All \(3\)-ring spaces \(V_3(q_i)\)
coincide with \(3\)-Selmer space \(V=V_3\)
\cite[Thm. 3.2, p. 2215]{Ma2014}.
\end{proof}


\subsection{Classification of Angell's \(3169\) simply real cubic fields}
\label{ss:Angell1972}

\noindent
In order to demonstrate the powerful performance
of our innovative techniques,
we construct all \(3\)-ring class fields \(K_{3,f}\)
which contain the normal closures \(N\) of the simply real cubic fields \(L\)
in Angell's table
\cite{Angell1972,Angell1973}
as abelian extensions of the associated imaginary quadratic base fields \(K<N\).

There arise four values of the \textit{multiplicity} \(m=1,2,3,4\),
and accordingly simply real cubic fields are collected in
singlets, doublets, triplets and quartets.
\textit{Nilets} with \(m=0\) complete the view.

The classification of the pure cubic fields,
respectively non-pure simply real cubic fields,
into \textbf{differential principal factorization types}
was established in
\cite{AMITA2020},
respectively
\cite{Ma2019b}.

Although the types \(\alpha\) and \(\beta\) of pure cubic fields
are similar to the types \(\alpha_2\) and \(\beta\) of non-pure simply real cubic fields,
we do not mix the classifications,
since firstly the existence of radicals among the principal factors
distinguishes pure cubic fields from non-pure simply real cubic fields,
and secondly, type \(\gamma\) can only occur for the former,
whereas type \(\alpha_1\) is only possible for the latter.


\renewcommand{\arraystretch}{1.1}

\begin{table}[ht]
\caption{Cubic discriminants in the range \(-2\cdot 10^4<d_L=f^2\cdot d<0\)}
\label{tbl:AngellSimplyReal}
\begin{center}
\begin{tabular}{|rl||r|rrrrr||rrr|}
\hline
       &           &       & \multicolumn{5}{c||}{Multiplicity}     & \multicolumn{3}{c|}{DPF} \\
 \(f\) & Condition & Total & \(0\) & \(1\) & \(2\) & \(3\) & \(4\) & \(\alpha_1\) & \(\alpha_2\) & \(\beta\) \\
\hline
 \(q\)      &   \(\equiv -1\,(\mathrm{mod}\,3)\) &  \(454\) &    \(\) &  \(454\) &  \(\) &   \(\) &   \(\) &     \(\) &   \(\) & \(454\) \\
 \(3\)      &  \(d\equiv +3\,(\mathrm{mod}\,9)\) &   \(62\) &    \(\) &   \(62\) &  \(\) &   \(\) &   \(\) &     \(\) &   \(\) &  \(62\) \\
 \(3\)      &  \(d\equiv -3\,(\mathrm{mod}\,9)\) &   \(58\) &    \(\) &   \(58\) &  \(\) &   \(\) &   \(\) &     \(\) &   \(\) &  \(58\) \\
 \(9\)      &  \(d\equiv -3\,(\mathrm{mod}\,9)\) &    \(7\) &    \(\) &     \(\) &  \(\) &  \(7\) &   \(\) &     \(\) &   \(\) &  \(21\) \\
 \(9\)      &  \(d\equiv -1\,(\mathrm{mod}\,3)\) &   \(23\) &    \(\) &   \(23\) &  \(\) &   \(\) &   \(\) &     \(\) &   \(\) &  \(23\) \\
\hline
 \(9\)      &  \(d\equiv +1\,(\mathrm{mod}\,3)\) &   \(20\) &    \(\) &   \(20\) &  \(\) &   \(\) &   \(\) &     \(\) & \(16\) &   \(4\) \\
 \(\ell\)   &   \(\equiv +1\,(\mathrm{mod}\,3)\) &   \(64\) &    \(\) &   \(64\) &  \(\) &   \(\) &   \(\) &     \(\) & \(49\) &  \(15\) \\
\hline
 \(q_1q_2\) &   \(\equiv -1\,(\mathrm{mod}\,3)\) &    \(6\) &    \(\) &     \(\) & \(6\) &   \(\) &   \(\) &     \(\) &   \(\) &  \(12\) \\
 \(3q\)     &  \(d\equiv +3\,(\mathrm{mod}\,9)\) &    \(7\) &    \(\) &     \(\) & \(7\) &   \(\) &   \(\) &     \(\) &   \(\) &  \(14\) \\
 \(3q\)     &  \(d\equiv -3\,(\mathrm{mod}\,9)\) &    \(3\) &    \(\) &     \(\) & \(3\) &   \(\) &   \(\) &     \(\) &   \(\) &   \(6\) \\
 \(9q\)     &  \(d\equiv -1\,(\mathrm{mod}\,3)\) &    \(3\) &    \(\) &     \(\) & \(3\) &   \(\) &   \(\) &     \(\) &   \(\) &   \(6\) \\
\hline
 \(9q\)     &  \(d\equiv +1\,(\mathrm{mod}\,3)\) &    \(3\) &    \(\) &     \(\) & \(3\) &   \(\) &   \(\) &     \(\) &   \(\) &   \(6\) \\
 \(3\ell\)  &  \(d\equiv +3\,(\mathrm{mod}\,9)\) &    \(1\) &    \(\) &     \(\) & \(1\) &   \(\) &   \(\) &     \(\) &   \(\) &   \(2\) \\
 \(q\ell\)  & \(\equiv\mp 1\,(\mathrm{mod}\,3)\) &    \(1\) &    \(\) &     \(\) & \(1\) &   \(\) &   \(\) &     \(\) &   \(\) &   \(2\) \\
\hline
\hline
 \(1\)      &                    \(\varrho_3=1\) & \(2143\) &    \(\) & \(2143\) &  \(\) &   \(\) &   \(\) & \(2143\) &   \(\) &    \(\) \\
\hline
 \(q\)      &   \(\equiv -1\,(\mathrm{mod}\,3)\) &  \(196\) & \(162\) &     \(\) &  \(\) & \(34\) &   \(\) &   \(87\) &   \(\) &  \(15\) \\
 \(3\)      &  \(d\equiv +3\,(\mathrm{mod}\,9)\) &   \(24\) &  \(22\) &     \(\) &  \(\) &  \(2\) &   \(\) &    \(4\) &   \(\) &   \(2\) \\
 \(3\)      &  \(d\equiv -3\,(\mathrm{mod}\,9)\) &   \(22\) &  \(16\) &     \(\) &  \(\) &  \(6\) &   \(\) &   \(13\) &   \(\) &   \(5\) \\
 \(9\)      &  \(d\equiv -1\,(\mathrm{mod}\,3)\) &    \(5\) &   \(5\) &     \(\) &  \(\) &   \(\) &   \(\) &     \(\) &   \(\) &    \(\) \\
\hline
 \(9\)      &  \(d\equiv +1\,(\mathrm{mod}\,3)\) &    \(9\) &   \(8\) &     \(\) &  \(\) &  \(1\) &   \(\) &    \(2\) &   \(\) &   \(1\) \\
 \(\ell\)   &   \(\equiv +1\,(\mathrm{mod}\,3)\) &   \(22\) &  \(19\) &     \(\) &  \(\) &  \(3\) &   \(\) &    \(7\) &   \(\) &   \(2\) \\
\hline
 \(q_1q_2\) &   \(\equiv -1\,(\mathrm{mod}\,3)\) &    \(2\) &   \(1\) &     \(\) &  \(\) &  \(1\) &   \(\) &     \(\) &   \(\) &   \(3\) \\
 \(3q\)     &  \(d\equiv +3\,(\mathrm{mod}\,9)\) &    \(3\) &   \(1\) &     \(\) &  \(\) &  \(2\) &   \(\) &     \(\) &   \(\) &   \(6\) \\
\hline
 \(9q\)     &  \(d\equiv +1\,(\mathrm{mod}\,3)\) &    \(1\) &    \(\) &     \(\) &  \(\) &  \(1\) &   \(\) &     \(\) &   \(\) &   \(3\) \\
 \(q\ell\)  & \(\equiv\mp 1\,(\mathrm{mod}\,3)\) &    \(2\) &   \(1\) &     \(\) &  \(\) &  \(1\) &   \(\) &     \(\) &   \(\) &   \(3\) \\
\hline
\hline
 \(1\)      &                    \(\varrho_3=2\) &   \(22\) &    \(\) &     \(\) &  \(\) &   \(\) & \(22\) &   \(88\) &   \(\) &    \(\) \\
\hline
 & Summary                                       & \(3163\) & \(235\) & \(2824\) &\(24\) & \(58\) & \(22\) & \(2344\) & \(65\) & \(725\) \\
\hline
\end{tabular}
\end{center}
\end{table}

\noindent
\textbf{Results:}
According to Table
\ref{tbl:AngellSimplyReal},
the number of all non-pure simply real cubic fields \(L\) having discriminants \(-2\cdot 10^4<d_L<0\) is given by \(\mathbf{3134}\).
Together with \(35\) pure cubic fields in Table
\ref{tbl:AngellPure},
the total number is \(\mathbf{3169}\), as announced correctly in
\cite{Angell1973}.

We emphasize the difference between
the \textit{number of discriminants} (without multiplicities)
\[2824+24+58+22=2928,\]
and the \textit{number of fields} (including multiplicities in a weighted sum)
\[1\cdot 2824+2\cdot 24+3\cdot 58+4\cdot 22=2824+48+174+88=3134,\]
which can be confirmed by adding the contributions to the \(3\) DPF types \(\alpha_1\), \(\alpha_2\), \(\beta\)
\[2344+65+725=3134.\]

In contrast, \(235\) is the number of \textit{formal cubic discriminants}
\(d_L=f^2\cdot d_K\) with fundamental discriminants \(d_K\) of imaginary quadratic fields
and \(3\)-admissible conductors \(f\) for each \(K\),
where the relevant multiplicity formula
\cite{Ma2014}
yields the value zero.
So the formal cubic discriminants belong to \textit{nilets},
i.e., multiplets with multiplicity \(m_3(K,f)=0\).
The total number of all (actual) cubic discriminants and formal cubic discriminants
is the number of admissible cubic discriminants
\[2928+235=3163.\]

According to Theorem
\ref{thm:ComplexRank0},
\textit{Nilets} can only arise for \(\varrho_3\ge 1\), but not for \(\varrho_3=0\).


\renewcommand{\arraystretch}{1.1}

\begin{table}[ht]
\caption{Pure cubic discriminants in the range \(-2\cdot 10^4<d_L=-3\cdot f^2<0\)}
\label{tbl:AngellPure}
\begin{center}
\begin{tabular}{|rl||r|rrrrr||rrr|}
\hline
       &           &       & \multicolumn{5}{c||}{Multiplicity}     & \multicolumn{3}{c|}{DPF} \\
 \(f\) & Condition & Total & \(0\) & \(1\) & \(2\) & \(3\) & \(4\) & \(\alpha\) & \(\beta\) & \(\gamma\) \\
\hline
 \(q\)          &   \(\equiv -1\,(\mathrm{mod}\,3)\) &  \(11\) &   \(8\) &  \(3\) &  \(\) &   \(\) &   \(\) &   \(\) &   \(\) & \(3\) \\
 \(9\)          &   \(d=-3\)                         &   \(1\) &    \(\) &  \(1\) &  \(\) &   \(\) &   \(\) &   \(\) &   \(\) & \(1\) \\
\hline
 \(\ell\)       &   \(\equiv +1\,(\mathrm{mod}\,3)\) &  \(10\) &   \(7\) &  \(3\) &  \(\) &   \(\) &   \(\) &  \(3\) &   \(\) &  \(\) \\
\hline
 \(q_1q_2\)     &   \(\equiv -1\,(\mathrm{mod}\,3)\) &   \(6\) &   \(1\) &  \(5\) &  \(\) &   \(\) &   \(\) &   \(\) &  \(5\) &  \(\) \\
 \(3q\)         &   \(d=-3\)                         &   \(5\) &   \(1\) &  \(4\) &  \(\) &   \(\) &   \(\) &   \(\) &  \(4\) &  \(\) \\
 \(9q\)         &   \(d=-3\)                         &   \(2\) &    \(\) &   \(\) & \(2\) &   \(\) &   \(\) &   \(\) &  \(4\) &  \(\) \\
\hline
 \(3\ell\)      &   \(d=-3\)                         &   \(3\) &   \(1\) &  \(2\) &  \(\) &   \(\) &   \(\) &  \(2\) &   \(\) &  \(\) \\
 \(9\ell\)      &   \(d=-3\)                         &   \(1\) &    \(\) &   \(\) & \(1\) &   \(\) &   \(\) &  \(2\) &   \(\) &  \(\) \\
 \(q\ell\)      & \(\equiv\mp 1\,(\mathrm{mod}\,3)\) &   \(8\) &   \(2\) &  \(6\) &  \(\) &   \(\) &   \(\) &  \(4\) &  \(2\) &  \(\) \\
\hline
 \(q_1q_2\ell\) & \(\equiv\mp 1\,(\mathrm{mod}\,3)\) &   \(1\) &    \(\) &  \(1\) &  \(\) &   \(\) &   \(\) &   \(\) &  \(1\) &  \(\) \\
 \(3q_1q_2\)    &   \(d=-3\)                         &   \(2\) &    \(\) &  \(2\) &  \(\) &   \(\) &   \(\) &   \(\) &  \(2\) &  \(\) \\
\hline
 \(3q\ell\)     &   \(d=-3\)                         &   \(2\) &    \(\) &  \(2\) &  \(\) &   \(\) &   \(\) &   \(\) &  \(2\) &  \(\) \\
\hline
 & Summary                                           &  \(52\) &  \(20\) & \(29\) & \(3\) &   \(\) &   \(\) & \(11\) & \(20\) & \(4\) \\
\hline
\end{tabular}
\end{center}
\end{table}

\noindent
According to Table
\ref{tbl:AngellPure},
the number of pure cubic fields \(L\) with discriminant \(-2\cdot 10^4<d_L<0\) is \(35\). \\
Actually, triplets and quartets of pure cubic fields do not occur in this range.

There is a difference between
the \textit{number of discriminants} (without multiplicities)
\[29+3=32,\]
and the \textit{number of fields} (including multiplicities in a weighted sum)
\[1\cdot 29+2\cdot 3=29+6=35,\]
which can be confirmed by adding the contributions to the \(3\) DPF types
\[11+20+4=35.\]

The total number of all (actual) cubic discriminants and formal cubic discriminants
(of the \(20\) nilets)
is the number of admissible pure cubic discriminants \(d_L=-3\cdot f^2\),
\[32+20=52.\]


\section{Multiplets over real quadratic fields for \(p=3\)}
\label{s:MultipletsReal3}

\noindent
We continue with real quadratic base fields \(K\)
having elevated \(3\)-Selmer rank \(\sigma_3=\varrho_3+1\),
due to the existence of a torsion free fundamental unit \(\varepsilon>1\).


\subsection{Classification of Angell's \(4804\) totally real cubic fields}
\label{ss:Angell1975}

\noindent
In order to demonstrate our progressive perspective
of classification of heterogeneous multiplets \(\mathbf{M}(K_{3,f})\)
into an enigmatic variety of differential principal factorization types,
we construct all \(3\)-ring class fields \(K_{3,f}\)
which contain the normal closures \(N\) of the totally real cubic fields \(L\)
in Angell's table
\cite{Angell1975,Angell1976}
as abelian extensions of the associated real quadratic base fields \(K<N\).

Again there arise four values of the \textit{multiplicity} \(m=1,2,3,4\),
and accordingly totally real cubic fields are collected in
singlets, doublets, triplets and quartets.
Formal \textit{nilets} complete the view.

The classification into \textbf{differential principal factorization types}
for non-cyclic totally real cubic fields was developed in
\cite{Ma1991b,Ma1991c,Ma2019b}.


\renewcommand{\arraystretch}{1.1}

\begin{table}[ht]
\caption{Cubic discriminants in the range \(0<d_L=f^2\cdot d<10^5\)}
\label{tbl:AngellTotallyReal}
\begin{center}
\begin{tabular}{|rl||r|rrrrr||rrrr|rrr|}
\hline
       &           &       & \multicolumn{5}{c||}{Multiplicity}     & \multicolumn{7}{c|}{Differential Principal Factorization} \\
 \(f\) & Condition & Total & \(0\) & \(1\) & \(2\) & \(3\) & \(4\) & \(\alpha_1\) & \(\beta_1\) & \(\beta_2\) & \(\gamma\) & \(\delta_1\) & \(\delta_2\) & \(\varepsilon\) \\
\hline
 \(q\)      &   \(\equiv -1\,(\mathrm{mod}\,3)\) & \(3025\) & \(2219\) & \(806\) &  \(\) &   \(\) & \(\) & \(\) &  \(\) &   \(\) &   \(\) &   \(\) &   \(\) & \(806\) \\
 \(3\)      &  \(d\equiv +3\,(\mathrm{mod}\,9)\) &  \(396\) &  \(287\) & \(109\) &  \(\) &   \(\) & \(\) & \(\) &  \(\) &   \(\) &   \(\) &   \(\) &   \(\) & \(109\) \\
 \(3\)      &  \(d\equiv -3\,(\mathrm{mod}\,9)\) &  \(389\) &  \(284\) & \(105\) &  \(\) &   \(\) & \(\) & \(\) &  \(\) &   \(\) &   \(\) &   \(\) &   \(\) & \(105\) \\
 \(9\)      &  \(d\equiv -3\,(\mathrm{mod}\,9)\) &   \(48\) &    \(9\) &  \(38\) &  \(\) &  \(1\) & \(\) & \(\) &  \(\) &   \(\) &   \(\) &   \(\) &   \(\) &  \(41\) \\
 \(9\)      &  \(d\equiv -1\,(\mathrm{mod}\,3)\) &  \(136\) &  \(102\) &  \(34\) &  \(\) &   \(\) & \(\) & \(\) &  \(\) &   \(\) &   \(\) &   \(\) &   \(\) &  \(34\) \\
\hline
 \(9\)      &  \(d\equiv +1\,(\mathrm{mod}\,3)\) &  \(127\) &   \(96\) &  \(31\) &  \(\) &   \(\) & \(\) & \(\) &  \(\) &  \(8\) &   \(\) &   \(\) & \(20\) &   \(3\) \\
 \(\ell\)   &   \(\equiv +1\,(\mathrm{mod}\,3)\) &  \(402\) &  \(316\) &  \(86\) &  \(\) &   \(\) & \(\) & \(\) &  \(\) & \(20\) &   \(\) &   \(\) & \(59\) &   \(7\) \\
\hline
 \(q_1q_2\) &   \(\equiv -1\,(\mathrm{mod}\,3)\) &   \(70\) &   \(30\) &  \(38\) & \(2\) &   \(\) & \(\) & \(\) &  \(\) &   \(\) & \(38\) &   \(\) &   \(\) &   \(4\) \\
 \(3q\)     &  \(d\equiv +3\,(\mathrm{mod}\,9)\) &   \(46\) &   \(23\) &  \(23\) &  \(\) &   \(\) & \(\) & \(\) &  \(\) &   \(\) & \(23\) &   \(\) &   \(\) &    \(\) \\
 \(3q\)     &  \(d\equiv -3\,(\mathrm{mod}\,9)\) &   \(45\) &   \(19\) &  \(25\) & \(1\) &   \(\) & \(\) & \(\) &  \(\) &   \(\) & \(25\) &   \(\) &   \(\) &   \(2\) \\
 \(9q\)     &  \(d\equiv -3\,(\mathrm{mod}\,9)\) &    \(5\) &     \(\) &    \(\) & \(4\) &  \(1\) & \(\) & \(\) &  \(\) &   \(\) &  \(9\) &   \(\) &   \(\) &   \(2\) \\
 \(9q\)     &  \(d\equiv -1\,(\mathrm{mod}\,3)\) &   \(14\) &    \(6\) &   \(8\) &  \(\) &   \(\) & \(\) & \(\) &  \(\) &   \(\) &  \(8\) &   \(\) &   \(\) &    \(\) \\
\hline
 \(9q\)     &  \(d\equiv +1\,(\mathrm{mod}\,3)\) &   \(15\) &    \(5\) &  \(10\) &  \(\) &   \(\) & \(\) & \(\) &  \(\) & \(10\) &   \(\) &   \(\) &   \(\) &    \(\) \\
 \(9\ell\)  &  \(d\equiv -1\,(\mathrm{mod}\,3)\) &    \(1\) &     \(\) &   \(1\) &  \(\) &   \(\) & \(\) & \(\) &  \(\) &  \(1\) &   \(\) &   \(\) &   \(\) &    \(\) \\
 \(3\ell\)  &  \(d\equiv +3\,(\mathrm{mod}\,9)\) &    \(6\) &    \(1\) &   \(5\) &  \(\) &   \(\) & \(\) & \(\) &  \(\) &  \(5\) &   \(\) &   \(\) &   \(\) &    \(\) \\
 \(3\ell\)  &  \(d\equiv -3\,(\mathrm{mod}\,9)\) &    \(5\) &    \(2\) &   \(3\) &  \(\) &   \(\) & \(\) & \(\) &  \(\) &  \(3\) &   \(\) &   \(\) &   \(\) &    \(\) \\
 \(q\ell\)  & \(\equiv\mp 1\,(\mathrm{mod}\,3)\) &   \(43\) &   \(13\) &  \(29\) & \(1\) &   \(\) & \(\) & \(\) &  \(\) & \(29\) &   \(\) &   \(\) &   \(\) &   \(2\) \\
\hline
\(3q_1q_2\) &   \(\equiv -1\,(\mathrm{mod}\,3)\) &    \(2\) &     \(\) &   \(1\) & \(1\) &   \(\) & \(\) & \(\) &  \(\) &   \(\) &  \(3\) &   \(\) &   \(\) &    \(\) \\
\hline
\hline
 \(1\)      &                    \(\varrho_3=1\) & \(3300\) &     \(\) &\(3300\) &  \(\) &   \(\) & \(\) & \(\) &  \(\) &   \(\) &   \(\) &\(3300\)&   \(\) &    \(\) \\
\hline
 \(q\)      &   \(\equiv -1\,(\mathrm{mod}\,3)\) &  \(275\) &  \(261\) &    \(\) &  \(\) & \(14\) & \(\) & \(\) & \(4\) &   \(\) &   \(\) & \(36\) &   \(\) &   \(2\) \\
 \(3\)      &  \(d\equiv -3\,(\mathrm{mod}\,9)\) &   \(35\) &   \(34\) &    \(\) &  \(\) &  \(1\) & \(\) & \(\) &  \(\) &   \(\) &   \(\) &  \(3\) &   \(\) &    \(\) \\
 \(\ell\)   &   \(\equiv +1\,(\mathrm{mod}\,3)\) &   \(28\) &   \(25\) &    \(\) &  \(\) &  \(3\) & \(\) & \(\) & \(3\) &   \(\) &   \(\) &  \(6\) &   \(\) &    \(\) \\
\hline
 \(3q\)     &  \(d\equiv -3\,(\mathrm{mod}\,9)\) &    \(2\) &    \(1\) &    \(\) &  \(\) &  \(1\) & \(\) & \(\) & \(3\) &   \(\) &   \(\) &   \(\) &   \(\) &    \(\) \\
\hline
\hline
 \(1\)      &                    \(\varrho_3=2\) &    \(5\) &     \(\) &    \(\) &  \(\) &   \(\) &\(5\) &\(16\)&  \(\) &   \(\) &   \(\) &  \(4\) &   \(\) &    \(\) \\
\hline
\hline
 & Summary                                       & \(8420\) & \(3733\) &\(4652\) & \(9\) & \(21\) &\(5\) &\(16\)&\(10\) & \(76\) &\(106\) &\(3349\)& \(79\) &\(1117\) \\
\hline
\end{tabular}
\end{center}
\end{table}

\noindent
\textbf{Results:}
According to Table
\ref{tbl:AngellTotallyReal},
the number of non-cyclic totally real cubic fields \(L\) with discriminant \(0<d_L<10^5\) is \(\mathbf{4753}\),
in perfect accordance with the results by Llorente and Oneto
\cite{LlOn1980,LlOn1982},
who discovered the ommission of ten fields in the table by Angell
\cite{Angell1975,Angell1976}.
Together with \(51\) cyclic cubic fields in Table
\ref{tbl:AngellCyclic},
the total number is \(\mathbf{4804}\)
(not \(4794\), as announced erroneously in
\cite{Angell1976}).

Again we emphasize the difference between
the \textit{number of discriminants} (without multiplicities)
\[4652+9+21+5=4687,\]
and the \textit{number of fields} (including multiplicities in a weighted sum)
\[1\cdot 4652+2\cdot 9+3\cdot 21+4\cdot 5=4652+18+63+20=4753,\]
which can be confirmed by adding the contributions to the \(7\) DPF types
(\(\alpha_2\), \(\alpha_3\) do not occur)
\[16+10+76+106+3349+79+1117=4753.\]
In contrast, \(3733\) is the number of \textit{formal cubic discriminants}
\(d_L=f^2\cdot d_K\) with fundamental discriminants \(d_K\) of real quadratic fields
and \(3\)-admissible conductors \(f\) for each \(K\),
where the relevant multiplicity formula
\cite{Ma2014}
yields the value zero.
So the formal cubic discriminants belong to \textit{nilets},
i.e., multiplets with multiplicity \(m_3(K,f)=0\).
The total number of all (actual) cubic discriminants and formal cubic discriminants
is the number of admissible cubic discriminants
\[4687+3733=8420.\]


\renewcommand{\arraystretch}{1.1}

\begin{table}[ht]
\caption{Cyclic cubic discriminants in the range \(0<d_L=f^2<10^5\)}
\label{tbl:AngellCyclic}
\begin{center}
\begin{tabular}{|rl||rr||r|}
\hline
                  &                          & \multicolumn{2}{c||}{M} & \multicolumn{1}{c|}{DPF} \\
 \(f\)            & Condition                        &  \(1\) &  \(2\) & \(\zeta\) \\
\hline
 \(9\)            & \(d=1\)                          &  \(1\) &   \(\) &  \(1\) \\
 \(\ell\)         & \(\equiv +1\,(\mathrm{mod}\,3)\) & \(30\) &   \(\) & \(30\) \\
\hline
 \(9\ell\)        & \(d=1\)                          &   \(\) &  \(4\) &  \(8\) \\
 \(\ell_1\ell_2\) & \(\equiv +1\,(\mathrm{mod}\,3)\) &   \(\) &  \(6\) & \(12\) \\
\hline
                  & Summary                          & \(31\) & \(10\) & \(51\) \\
\hline
\end{tabular}
\end{center}
\end{table}

\noindent
According to Table
\ref{tbl:AngellCyclic},
the number of cyclic cubic fields \(L\) with discriminant \(0<d_L<10^5\) is \(51\),
with \(31\) arising from singlets having conductors \(f\) with a single prime divisor,
and \(20\) from doublets having two prime divisors of the conductor \(f\).
(M denotes the multiplicity.)

We point out that cyclic cubic fields
are rather contained in \textit{ray class fields} over \(\mathbb{Q}\)
than in ring class fields over real quadratic base fields.
The single possible DPF type \(\zeta\)
has nothing to do with the \(9\) DPF types
\(\alpha_1,\alpha_2,\alpha_3,\beta_1,\beta_2,\gamma,\delta_1,\delta_2,\varepsilon\)
of non-abelian totally real cubic fields in
\cite{Ma2019b}.


\section{Conclusion and outlook}
\label{s:Conclusion}

\noindent
In this paper, we have classified all multiplets
\(\mathrm{Obj}(K_{3,f})\)
of \textit{non-pure simply real} cubic fields \(L\)
(more precisely of their normal closures \(N\))
according to the associated multiplets of invariants,
namely the \textbf{differential principal factorization types},
\(\mathrm{Inv}(K_{3,f})\),
where \(K_{3,f}\) denotes the \(3\)-ring class field
modulo a \(3\)-admissible conductor \(f\) of the imaginary quadratic subfield \(K<N\): \\
(Recall that \(\mathrm{Obj}(K_{3,f})=(N_{f,i})_{1\le i\le m}\)
and \(\mathrm{Inv}(K_{3,f})=(\tau(N_{f,i}))_{1\le i\le m}\),
here \textit{homogeneously}.) \\
\(2824\) \textit{singulets} of type either \((\alpha_1)\) or \((\alpha_2)\) or \((\beta)\), according to Table
\ref{tbl:AngellSimplyReal}; \\
\(24\) \textit{doublets} of exclusive type \((\beta,\beta)\) (without \(3\) pure cubic doublets); \\
\(58\) \textit{triplets} with the following distribution of types:

\(7\) triplets of type \((\beta,\beta,\beta)\) for \(f=9\) singular, \(\varrho_3=0\),

and \(51\) triplets sharing common \(3\)-class rank \(\varrho_3=1\) of \(K\)
\cite[Tbl. 1, pp. 118--121]{Sm1982},
namely

\(34\) triplets of type \((\alpha_1,\alpha_1,\alpha_1)\) for \(f=q,\ell,3\),

\(3\) triplets of type \((\alpha_1,\alpha_1,\beta)\) for \(f=\ell,9\) split,

\(5\) triplets of type \((\alpha_1,\beta,\beta)\) for \(f=q,3\), and

\(9\) triplets of type \((\beta,\beta,\beta)\) for \(f=q,3,3q,q\ell,9q\), finer than
\cite{Sm1982} since \(\alpha_2\) does not occur; and \\
\(22\) \textit{quartets} of exclusive type \((\alpha_1,\alpha_1,\alpha_1,\alpha_1)\)
(see
\cite{Ma1991a}
for details concerning the capitulation).

\bigskip
Similarly, we have classified all multiplets
\(\mathrm{Obj}(K_{3,f})\)
of \textit{non-cyclic totally real} cubic fields \(L\)
(more precisely of their normal closures \(N\))
according to the associated multiplets of invariants,
namely the \textbf{differential principal factorization types},
\(\mathrm{Inv}(K_{3,f})\),
where \(K_{3,f}\) denotes the \(3\)-ring class field
modulo a \(3\)-admissible conductor \(f\) of the real quadratic subfield \(K<N\): \\
\(4652\) \textit{singulets} of type either \((\beta_2)\) or \((\gamma)\) or \((\delta_1)\) or \((\delta_2)\) or \((\varepsilon)\), according to Table
\ref{tbl:AngellTotallyReal}; \\
\(9\) \textit{doublets},
\(4\) of type \((\gamma,\gamma)\) for \(f=9q,3q_1q_2\) and
\(5\) of type \((\varepsilon,\varepsilon)\) for \(f=3q,9q,q_1q_2,q\ell\); \\
\(21\) \textit{triplets} with the following distribution of types:

\(1\) triplet of type \((\varepsilon,\varepsilon,\varepsilon)\) for \(f=9\) singular, \(\varrho_3=0\),

\(1\) triplet of type \((\gamma,\gamma,\gamma)\) for \(f=9q\) singular, \(\varrho_3=0\),

and \(19\) triplets sharing common \(3\)-class rank \(\varrho_3=1\) of \(K\)
(with considerable refinement of

Schmithals' coarse distinction of only two alternatives
\cite[Tbl. 2, pp. 122--123]{Sm1982}),
namely

\(13\) triplets of type \((\delta_1,\delta_1,\delta_1)\) for \(f=3,q\),

\(1\) triplet of type \((\beta_1,\beta_1,\beta_1)\) for \(f=3q\),

\(2\) triplets of type \((\beta_1,\beta_1,\varepsilon)\) for \(f=q\)
(conspicuously with symbol \lq\lq\(-\)\rq\rq\ in
\cite{Sm1982}),
and

\(3\) triplets of type \((\beta_1,\delta_1,\delta_1)\) for \(f=\ell\), \\
\(5\) \textit{quartets},
\(1\) of type \((\alpha_1,\alpha_1,\alpha_1,\alpha_1)\) and
\(4\) of type \((\alpha_1,\alpha_1,\alpha_1,\delta_1)\)
(more details in
\cite{Ma2012,Ma2014b}).

\bigskip
In the same manner, we shall \textit{refine more extensive tables} by
Fung and Williams
\cite{FuWi1990},
Ennola and Turunen
\cite{EnTu1983,EnTu1985},
Llorente and Quer
\cite{LlQu1988}
in the new year \(2021\).

Moreover, we shall provide extensive evidence
of the \textit{truth of Scholz' conjecture},
which we have proved for \(p=3\) in
\cite{Ma2019b},
also for \(p=5\) and \(p=7\),
and probably for any odd prime \(p\).

Implementations of our innovative algorithms in Magma
\cite{BCP1997,BCFS2020,MAGMA2020}
may be requested via email.

\newpage

\section{Acknowledgements}
\label{s:Thanks}

\noindent
This work has been completed on Tuesday, 29 December 2020.
In order \textit{to disprove any claims of priority}
concerning the innovative perspective of
classifying \textit{multiplets of dihedral fields},
contained in \(p\)-ring class fields,
into \textit{differential principal factorization types},
the article has immediately been disseminated
on various scientific open access platforms.

\noindent
The author gratefully acknowledges
that his research was supported by the Austrian Science Fund (FWF):
projects J0497-PHY and P26008-N25.



\end{document}